\documentclass[11pt]{amsart}

\usepackage{DKstyle}
\usepackage{bm}
\usepackage{microtype}
\usepackage{float}
\usepackage[toc,page]{appendix}

\newcommand{\vertiii}[1]{{\left\vert\kern-0.25ex\left\vert\kern-0.25ex\left\vert #1
    \right\vert\kern-0.25ex\right\vert\kern-0.25ex\right\vert}}
 \makeatletter
\newcommand*{\rom}[1]{\expandafter\@slowromancap\romannumeral #1@}
\makeatother
\theoremstyle{plain}

\usepackage{caption}
\usepackage[labelfont=rm]{subcaption}

\usepackage[pagebackref=false, colorlinks]{hyperref}

\hypersetup{pdffitwindow=true,linkcolor=blue,citecolor=blue,urlcolor=blue,menucolor=blue}

\usepackage{comment}

\def\d{\text{d}}

\subjclass{}

\date{\today}
\dedicatory{}
\commby{}
\thanks{S.Y was supported by the National Key R\&D Program of China (2024YFA1013603, 2024YFA1013600). D.K was supported by Simons grant SFI-MPS-TSM-00013031. }
\title{Second moment of incomplete  Eisenstein series revisited }
\author{Dubi Kelmer}
\address{Department of Mathematics, Boston College, Chestnut Hill MA 02467-3806, USA}
\email{kelmer@bc.edu}

\author{Shucheng Yu}
\address{School of Mathematical Sciences, University of Science and Technology of China (USTC), 230026, Hefei, China}
\email{yusc@ustc.edu.cn}

\begin{document}
\begin{abstract}
We revisit Yu's second moment formula for incomplete Eisenstein series and give uniform bounds for secondary terms appearing in this formula. We use these bounds to improve results on distribution of integer points on the light cone corresponding to an indefinite quadratic form,
as well as intrinsic Diophantine approximation on ellipsoids, and counting holonomy vectors of saddle connections of Veech surfaces.
\end{abstract}
\maketitle

\section{Introduction}
Fix an integer $n\geq1$. Let $Q$ be a real quadratic form of signature $(n+1,1)$ and let $G=\SO_Q^+(\R)$.
Fix an Iwasawa decomposition $G={U}AK$, where ${U}$ is unipotent, $A$ is a one-dimensional diagonalizable group, and $K$ is a maximal compact subgroup.
Let $M=Z_K(A)$ denote the centralizer of $A$ in $K$, let $L={U}M$ and let $\cV=L\bk G$.
Let $\mu_G$ and $\mu_L$ denote Haar measures on $G$ and $L$ respectively (normalized as in \eqref{equ:Haar1} and \eqref{equ:Lhaar} below) and let $m_{\cV}$ be the corresponding measure on $\cV=L\bk G$ normalized so that
$$
\d\mu_G=\d\mu_L \d m_{\cV}.
$$
Let  $\cV_Q^+$ denote one of the two connected components of the light cone defined by nonzero solutions to $Q(\bm{v})=0$. Let $\bm{e}_0\in \cV_Q^+$ be a point stabilized by $L$. The transitive action of $G$ allows us to identify $\cV_Q^+$ with $\cV=L\bk G$.  This identification gives a measure $m_{\cV_Q^+}=m_{\cV}$ on the light cone.

Let $\Gamma\leq G$ be a nonuniform lattice. After perhaps conjugating $\Gamma$ we may assume  that the parabolic group $P=UAM$ corresponds to a cusp of $\Gamma$, that is, that $\Gamma_L:=\Gamma\cap L$ is a lattice in $L$, and also that $\mu_L(\G_L\bk L)=1$.
{Throughout the paper, all functions are assumed to be measurable.} To any bounded, compactly supported function $f:\cV\to \C$, the corresponding incomplete Eisenstein series at $P$, denoted by
$E_f\in L^2(\Gamma\bk G)$, is given by
$$
E_f(g)=\sum_{\gamma \in \Gamma_L\bk \Gamma} f( \gamma g).
$$

Extending the results of \cite{KelmerMohammadi12} for $n=1,2$, Yu \cite{Yu17,YuThesis18} obtained moment formulas for these incomplete Eisenstein series. These formulas express the second moment through the space average of $f$ (with respect to $m_{\cV}$), a bounded-operator term, and secondary terms associated with exceptional poles. In the absence of exceptional poles, the second moment is bounded in terms of $\|f\|_1$ and $\|f\|_2$.  The secondary terms are multiples of $M_{f,f}(s_k)$, defined in \eqref{e:Ms}, with $s_k\in(\frac{n}{2},n)$ an exceptional pole of $\Gamma$.
Here for any $1\leq p\leq \infty$, $\|f\|_p$ denotes the $L^p$-norm of $f$ on $\cV$ with respect to $m_{\cV}$.

For many applications it is useful to have bounds for these secondary terms when $f=\chi_B$ is the indicator function of a measurable set $B\subseteq \cV$ in terms of the measure $m_\cV(B)$.
Recalling the identification of $\cV$ with the light cone $\cV_Q^+$, we note that when $B=B_t=tB_1$ is a dilation of a fixed relatively compact set $B_1\subseteq \cV_Q^+$ of positive measure, it is not hard to see that
$$
M_{\chi_{B_t},\chi_{B_t}}(s)=c m_{\cV_Q^+}(B_t)^{\frac{2s}{n}}
$$
for some constant $c>0$ depending on $s,n$ and $B_1$.
It is natural to expect that such a bound holds in general. While it was shown in \cite{YuThesis18,KelmerYu2023a} that this is indeed the case for a large family of sets (including  generalized sectors) we were not able to show this in complete generality, which was an obstacle  for some of these applications. Our main result here is the following bound for the secondary term $M_{f,f'}(s)$ for general bounded compactly supported functions on $\cV$.

\begin{Thm}\label{t:Msbound}
Let $\frac{n}{2}<s<n$. For any two bounded, compactly supported
functions $f,f':\cV\to \C$, we have
\begin{equation}\label{e:Msfbound}
|M_{f,f'}(s)|\ll_{n,s} \|f\|_\infty^{1-\frac{s}{n}}\|f\|_1^{\frac{s}{n}} \|f'\|_\infty^{1-\frac{s}{n}}\|f'\|_1^{\frac{s}{n}}.
\end{equation}
\end{Thm}

Here and below we write $A\ll B$ also $A=O(B)$ to indicate that $|A|\leq c |B|$ for some constant $c>0$, and we use subscripts to indicate that the constant $c$ may depend on additional parameters.
\begin{remark}
For $\lambda>0$, consider the dilations $f_\lambda(\bm v)=f(\lambda^{-1}\bm v)$ of a fixed nonzero nonnegative function $f$. It follows easily from the definition of $M_{f,f'}(s)$ that for any $\frac{n}{2}<s<n$ and for any $c,\lambda>0$,
\[
M_{cf_\lambda,cf_\lambda}(s)=c^2\lambda^{2s}M_{f,f}(s).
\]
This shows that the exponents in \thmref{t:Msbound} are best possible.
\end{remark}

With this bound we are able to extend the applications of \cite{KelmerYu2023a} to the distribution of points in $\Gamma$-orbits, to intrinsic Diophantine approximations on ellipsoids, and the application of \cite{BurrinFairchild2026} to counting  holonomy vectors of saddle connections of Veech surfaces.

\subsection{Application to counting integer points on the light cone}
\label{sec:recall}
Assume now that the quadratic form $Q$ is a rational and $\Q$-isotropic  form  of signature $(n+1,1)$. Denote by $\cV_Q^+$ one of the two connected components of the light cone and recall the identification of $\cV_Q^+$ with $\cV=L\bk G$.
Let $\cL_Q=\cV^+_Q\cap \Z^{n+2}_{\rm pr}$ denote the subset of primitive integer points in $\cV_Q^+$.
Let $\Gamma_Q=\SO_Q(\Z)\cap G$ and note that $\G_Q$ is a nonuniform lattice and that it acts on $\cL_Q$, splitting it into finitely many orbits.

For a bounded compactly supported function $f:\cV_Q^+\to\C$, its light-cone Siegel transform is the function on $\Gamma_Q\backslash G$ defined by
$$
S_Qf(g)=\sum_{\bm v\in\cL_Q}f(\bm vg).
$$
We use the same definition for nonnegative functions, allowing the value $+\infty$.
When taking $f(\bm{v})={\|\bm{e}_0\|_Q^s}\|\bm{v}\|_Q^{-s}$ with $\Re(s)>n$ and $\|\cdot\|_Q$  {the norm in \cite[(2.33)]{KelmerYu2023a}, obtained from the Euclidean norm by a linear change of variables},
we get the light-cone Eisenstein series
\begin{equation}\label{e:light-cone-eisenstein}
E_Q(s,g)=\|\bm e_0\|_Q^s\sum_{\bm{v}\in \cL_Q}\|\bm{v}g\|_Q^{-s}.
\end{equation}
The series absolutely converges when $\Re(s)>n$ and extends meromorphically to $\C$. In the half plane $\Re(s)\geq \frac{n}{2}$, it has a simple pole at $s=n$ with constant residue $\omega_Q>0$ and at most one additional exceptional pole at $s_n=\left\lfloor\frac{n+2}{2}\right\rfloor$ (when $n=1,2$ there are no exceptional poles).
{Write $m_{\cV_Q^+}(f)=\int_{\cV_Q^+}f\,\d m_{\cV_Q^+}$ and let $\mu_Q$ be the probability $G$-invariant measure on $\G_Q\bk G$. In \cite[Theorem 1.1]{KelmerYu2023a}, the first formula below is proved for a bounded compactly supported function $f$, and the second is proved under the additional assumption that $f$ is smooth.
}
\begin{equation}\label{e: momen1}
\int_{\G_Q\bk G} S_Qf(g)\,\d\mu_{Q}(g)=\omega_Qm_{\cV_Q^+}(f),
\end{equation}
and
\begin{equation}\label{e:moment1}
\int_{\G_Q\bk G} |S_Qf(g)|^2\,\d\mu_{Q}(g)=|\omega_Qm_{\cV_Q^+}(f)|^2+c_QM_{f,f}(s_n)+O(m_{\cV_Q^+}(|f|^2)).
\end{equation}
{Here when $n=1,2$, the Eisenstein series $E_Q(s,g)$ has no exceptional poles, and \eqref{e:moment1} holds without the secondary term. We set $c_Q=0$ in these two cases. When $n\geq3$, we have $\frac{n}{2}<s_n<n$ and $c_Q\geq 0$. The secondary term occurs precisely when $E_Q(s,g)$ has a pole at $s=s_n$, equivalently when $c_Q>0$.}
We also define the following exponent
\begin{align*}
\beta_Q:=\begin{cases}
1 &c_Q=0,
\\
2s_n/n &c_Q>0.
\end{cases}
\end{align*}

For a finite-measure Borel set $B\subseteq\cV_Q^+$, define the discrepancy function
$$
D(\cL_Q g,B)=|\#(\cL_Q g\cap B)-\omega_Q m_{\cV_Q^+}(B)|.
$$
{Combining \eqref{e:moment1} with \thmref{t:Msbound}, we obtain the following second moment estimate. Its discrepancy consequence extends \cite[Corollary 1.2]{KelmerYu2023a}.}

\begin{Thm}\label{t:Dis}
For any bounded, compactly supported function $f:\cV_Q^+\to\C$, we have
\begin{equation}\label{e:momentbound}
\int_{\Gamma_Q\bk G}|S_Qf(g)|^2\,\mathrm{d}\mu_Q(g)
=|\omega_Qm_{\cV_Q^+}(f)|^2
+O_Q\!\left(\|f\|_2^2+\|f\|_\infty^{2-\beta_Q}\|f\|_1^{\beta_Q}\right).
\end{equation}
If $c_Q=0$, the error term can be replaced by $O_Q(\|f\|_2^2)$.

In particular, there is a constant $c>0$ depending only on $Q$ such that, for every finite-measure Borel set $B\subseteq\cV_Q^+$,
\begin{equation}\label{e:DisMoment}
\int_{\Gamma_Q\bk G}D(\cL_Qg,B)^2\,\mathrm{d}\mu_Q(g)
\leq c\bigl(m_{\cV_Q^+}(B)+m_{\cV_Q^+}(B)^{\beta_Q}\bigr).
\end{equation}
\end{Thm}

\begin{remark}
The discrepancy bound \eqref{e:DisMoment} was previously proved for generalized sectors (see \cite[Definition 2.37]{KelmerYu2023a}) in \cite[Corollary 1.2]{KelmerYu2023a}. Here \eqref{e:DisMoment} extends it to every finite-measure Borel set, with a constant depending only on $Q$.
\end{remark}

As a direct consequence, we get the following bound that holds for almost all translates of a set, extending the results of  \cite[Theorem 1.3]{KelmerYu2023a} (see also \cite[Theorem 5]{YuThesis18} ) by removing the restriction of the allowable sets in cases when $c_Q\neq 0$.
\begin{Thm}\label{t:Dis2}
Let $\cB$ be an increasing family of finite-measure Borel sets in $\cV_Q^+$.
Then for any $\e>0$, for $\mu_Q$-a.e. $\G_Q g\in  \G_Q\backslash G$,  for all $B\in \cB$ with ${m_{\cV_Q^+}}(B)\geq 1$,
$$\#(\cL_Q g \cap B)=\omega_Q {m_{\cV_Q^+}}(B) +O_{\epsilon,\cB, g}\bigg({m_{\cV_Q^+}}(B)^{\frac{\beta_Q}{2}+\epsilon}\bigg).$$
\end{Thm}

\begin{remark}\label{r:lattice-surfaces}
We note that the action of $\Gamma_Q$ on $\cL_Q$ separates it into finitely many orbits, $\cL_1,\ldots,\cL_\kappa$ \cite[Lemma 2.2]{KelmerYu2023b}, corresponding to the cusps of $\Gamma_Q$, and we can obtain a similar result when counting in each orbit separately. More generally, given any nonuniform lattice $\Gamma\leq G$ with a cusp corresponding to a parabolic $P={U}AM$, that is, with $\Gamma_L=\Gamma\cap L$ a lattice in $L$, we can consider the $\Gamma$-orbit $\cL=L\Gamma$ in $\cV=L\bk G$, and get a similar result for counting the number of points in $\cL g\cap B$. The inner product estimate \eqref{e:incomplete-momentbound} for incomplete Eisenstein series below gives  an analogous result for $\#(\cL g\cap B)$.  The only difference is that in these cases the exponent $\beta_Q$ will be replaced by $\beta_\Gamma\in[1,2)$, {equal to $1$ if the corresponding Eisenstein series has no exceptional pole} and to $2s_\Gamma/n$ if its largest exceptional pole is $s_\Gamma$, {and the constant $\omega_Q$ is replaced by the explicit constant $\mu_L(\Gamma_L\backslash L)/\mu_G(\Gamma\backslash G)$} (see \thmref{t:Dis3} below).
\end{remark}

Given a function $F:\R^{n+2}\to \R^m$ one expects that, generically, the values $F(\bm{v})$ with $\bm{v}\in \cL_Q$ will be evenly distributed. This was shown to hold quantitatively  for a large family of polynomial maps   in  \cite[Theorem 4.3 and Theorem 4.4]{KelmerYu2023a} under the assumption that $c_Q=0$. As a consequence of the variance bound in \thmref{t:Dis}, we can now obtain the same results in general without this assumption.
For example in the special case when $F$ is linear we get the following.

\begin{Thm}\label{thm:linearformintro}
Assume $n\geq 2$. Let  $m<n$ be a positive integer and let $\{\Omega_T\}_{T>1}\subseteq \R^m$ be a decreasing family of bounded measurable sets with $\vol(\Omega_T)=T^{-a}$ for some $0<a<n-m$. Then there is some $\nu>0$ such that for almost every linear map $F:\R^{n+2}\to \R^m$ with $Q|_{\ker{F}}$ indefinite, for all sufficiently large $T$
$$
\#\{\bm{v}\in \cL_Q: \|\bm{v}\|_Q\leq T,\; F(\bm{v})\in \Omega_T\}=c(Q,F) T^{n-m-a}(1+O_F(T^{-\nu})),
$$
where $c(Q,F)$ is some positive constant depending only on $Q$ and $F$.
\end{Thm}
By taking the sets $\Omega_T$ to be shrinking norm balls in $\R^m$ we get the following corollary.
\begin{Cor}\label{c:GGN}
Assume $n\geq 2$. Let $m<n$ be a positive integer and  let $\kappa>\frac{m}{n-m}$. Given any $\bm\xi\in \R^m$, for almost every linear map $F:\R^{n+2}\to \R^m$ with $Q|_{\ker{F}}$ indefinite, for all $\epsilon>0$ sufficiently small,
the system of inequalities
\begin{equation}\label{e:GGN}
\|F(\bm{v}) -\bm\xi\|< \epsilon \quad \text{and}\quad \|\bm{v}\|< \epsilon^{-\kappa}
\end{equation}
has a solution $\bm{v}\in \cL_Q$.
\end{Cor}
\begin{remark}
These results extend \cite[Theorem 1.4 and Corollary 1.5]{KelmerYu2023a} by removing the assumption $c_Q=0$. {When $n=m+1$, the threshold $\kappa>m$ agrees with \cite[Theorem 1.3]{GhoshGorodnikNevo2020}. For $n>m+1$, the threshold $\kappa>\frac{m}{n-m}$ improves on that obtained by their method; see \cite[Remark 1.11]{KelmerYu2023a}}.
\end{remark}

\subsection{Application to intrinsic Diophantine approximation}

Let $Q_0$ be a positive definite rational quadratic form in $n+1$ variables and assume that $Q(\bm x,y)=Q_0(\bm x)-y^2$ is $\Q$-isotropic.
Let $\omega_Q,c_Q,\beta_Q$ be the parameters associated to this form as in Section \ref{sec:recall}. Write $\|\bm x\|_{Q_0}:=\sqrt{Q_0(\bm x)}$ and let
$$
S_{Q_0}:=\{\bm x\in\R^{n+1}:Q_0(\bm x)=1\}
$$
be the ellipsoid associated to $Q_0$. Let $\sigma_{Q_0}$ denote the $\SO_{Q_0}(\R)$ invariant probability surface measure on $S_{Q_0}$ (as defined in \cite[Section 1.3]{KelmerYu2023a}).

Given a function $\psi:\N\to (0,\infty)$ and $\bm{x}\in S_{Q_0}$, consider the counting function
\begin{align}\label{equ:khincountfunell}
\cN_{Q_0, \psi}(\bm{x}, T):=\#\left\{\tfrac{\bm{p}}{q}\in S_{Q_0}: \gcd(\bm{p},q)=1,\  \|\bm{x}-\tfrac{\bm{p}}{q}\|_{Q_0}<\psi(q),\ 1\leq q< T\right\}.
\end{align}

{Following the terminology of  \cite{KleinbockMerrill2015}, we call $\bm x\in S_{Q_0}$ $\psi$-approximable if $\|\bm x-\bm p/q\|_{Q_0}<\psi(q)$ holds for infinitely many distinct reduced rational points $\bm p/q\in S_{Q_0}$.} Note that the function $\cN_{Q_0,\psi}(\bm{x}, T)$ is increasing in $T$ and that it is unbounded {as $T\to\infty$} if and only if  $\bm{x}$ is $\psi$-approximable in $S_{Q_0}$.
We also write
\begin{align}\label{equ:partialsum}
 J_{\psi}(T):=\sum_{1\leq q< T}q^{n-1}\psi(q)^n\quad \mbox{and }\quad I_{\psi}(T):=\sum_{1\leq q< T}q^{n-1}\psi(q)^{n+2}.
 \end{align}
 We then have the following {quantitative Khintchine-type theorem removing the assumption that $c_Q=0$ in \cite[Theorem 1.9]{KelmerYu2023a}}.
\begin{Thm}\label{thm:khinttypegen}
{Let $Q(\bm{x},y)=Q_0(\bm{x})-y^2$ with $Q_0$ and $Q$ as in this subsection.} For any decreasing $\psi: \N\to (0,\infty)$ satisfying that $\lim\limits_{q\to\infty}\psi(q)=0$ and $\sum_{q\in\N}q^{n-1}\psi(q)^n=\infty$,  for $\sigma_{Q_0}$-a.e. $\bm{x}\in S_{Q_0}$ and for all sufficiently large $T$,
\begin{align*}
\cN_{Q_0,\psi}(\bm{x}, T)= {\frac{\omega_Q\Gamma((n+3)/2)}{\sqrt{\pi}(n+1)\Gamma((n+2)/2)}}J_{\psi}(T)+O_{Q_0, \psi}\left(J_{\psi}(T)^{\frac{n+2+\beta_Q}{n+4}}\log(J_{\psi}(T))+I_{\psi}(T)\right).
\end{align*}
\end{Thm}

As a consequence we get the following corollary extending the result of \cite[Theorem 6.3]{FishmanKleinbockMerrillSimmons2021}.

\begin{Cor}\label{c:1}
{Let $Q(\bm x,y)=Q_0(\bm x)-y^2$ with $Q_0$ and $Q$ as in this subsection.} For any decreasing function $\psi:\N\to (0,\infty)$ we have that $\sigma_{Q_0}$-almost every (respectively almost no) $\bm{x}\in S_{Q_0}$ is $\psi$-approximable if and only if the series $\sum_{q\in\N}q^{n-1}\psi(q)^n$ diverges (respectively converges).
\end{Cor}

\begin{remark}
\thmref{thm:khinttypegen} extends  \cite[Theorem 1.9]{KelmerYu2023a} by removing the assumption that $c_Q= 0$ (we note that when $c_Q=0$ then $\beta_Q=1$ and we recover \cite[Theorem 1.9]{KelmerYu2023a}). Also note that  \corref{c:1} extends the qualitative result of  \cite[Theorem 6.3]{FishmanKleinbockMerrillSimmons2021} by removing the assumption that the function $q\mapsto q\psi(q)$ is decreasing  (this was previously only proved when $c_Q=0$).
\end{remark}

\subsection{Application to counting saddle connections}
Let $M$ denote a Veech surface with $\Gamma_M\leq\SL_2(\R)$ its Veech group and let $S_M\subseteq \R^2$  denote the set of distinct holonomy vectors of its oriented saddle connections.
Let $\sigma_M\in [1/2,1)$ denote the largest exceptional pole of the Eisenstein series for $\Gamma_M$ (and set $\sigma_M=\tfrac{1}{2}$ if there are none).
In \cite[Theorem 4.1]{BurrinNevoRuhrWeiss2020} it was shown that for growing Euclidean balls $B_t=\{\bm{v}: \|\bm{v}\|\leq t\}$ for any $g\in \SL_2(\R)$
$$\#(S_M g\cap B_t)=c_M|B_t|+O(|B_t|^{\max\{\sigma_M,2/3\}}),$$
with  $c_M>0$ the corresponding Siegel--Veech constant and  $|B|$ denoting the planar Lebesgue measure of $B$.
In \cite[Theorem 10.3]{BurrinFairchild2026}, a different moment formula was used to estimate $\# (S_M g\cap B)$ for $B\subseteq \R^2$ in an arbitrary growing family of Borel sets, giving essentially the same bound for almost all $g\in {\GL_2(\R)}$.
By identifying $\SO^+(2,1)$ with  $\PSL_2(\R)$ and applying the inner product estimate to corresponding incomplete Eisenstein series (see \thmref{t:Dis3} below) we get the following.
\begin{Thm}\label{t:Veech}
Let $\cB$ be an increasing family of finite-measure Borel sets in $\R^2$.
Then  for any $\delta>\sigma_M$, for a.e. $g\in\SL_2(\R)$,  for all $B\in \cB$ with $|B|$ sufficiently large
\begin{align}\label{equ:sl2ct}
\#(S_M g \cap B)=c_M |B|+O_{\delta,M, \cB,g}\left(|B|^{\delta}\right).
\end{align}
\end{Thm}
\begin{rem}
By scalar dilation and Fubini's theorem, and by applying the same argument to the reflected family $\{BJ:B\in\cB\}$ with $J=\operatorname{diag}(-1,1)$
for the negative-determinant component,
\thmref{t:Veech} also gives, for any $\delta>\sigma_M$, for a.e.
$g\in\operatorname{GL}_2(\R)$ and for all $B\in\cB$ with $|B|\gg 1$,
\begin{align}\label{equ:gl2ct}
\#(S_Mg\cap B)|\det g|
=c_M|B|+O_{\delta,M,\cB,g}\bigl(|B|^\delta\bigr).
\end{align}
Taking transposes converts our row-vector convention to the
column-vector convention of \cite{BurrinFairchild2026}.
Their Theorem 10.3 gives the corresponding estimate for every
$\delta>\max\{2/3,\sigma_M\}$, so \eqref{equ:gl2ct} improves
the power-saving exponent when $\sigma_M<2/3$.
We note that this improvement follows from applying the counting argument to our variance bound $O_M(|B|+|B|^{2\sigma_M})$ on the corresponding quotient of $\PSL_2(\R)$, obtained from \eqref{e:incomplete-momentbound}, in place of the bound $O_M(|B|^{\max\{4/3,2\sigma_M\}})$ averaged over a cone in $\operatorname{GL}_2^+(\R)$ in \cite[Proposition 10.1]{BurrinFairchild2026}.
\end{rem}

\vspace{5pt}

\noindent\textbf{AI disclosure.}
ChatGPT (5.6 Sol and Astra) and Claude (Opus 5.5)  were used to assist with the preparation of this note. The authors reviewed the resulting material and take full responsibility for the content.

 \section{Preliminaries and notation}
 \subsection{Coordinates and measures}
Let $Q$ be a real quadratic form of signature $(n+1,1)$, and let $G=\SO_Q^+(\R)$.
Let  $\cV_Q^+$ denote one of the two connected components of the light cone defined by nonzero solutions to $Q(\bm{v})=0$.
We take the base point $\bm{e}_0\in\cV_Q^+$ and the corresponding Iwasawa decomposition $G={U}AK$ as in \cite[Section 2.1.1]{KelmerYu2023a}. We use the same parameters $u_{\bm x}\in {U}$, $\bm x\in\R^n$, and $a_y\in A$, $y>0$; in particular, $\bm e_0a_y=y^{-1}\bm e_0$.

Let $M=Z_K(A)\cong\SO(n)$, where $K\cong\SO(n+1)$. Then $L={U}M$ is the stabilizer of $\bm e_0$, and $P={U}AM$ is the stabilizer of the line $\R\bm e_0$. The map $Lg\mapsto\bm e_0g$ identifies $\cV=L\bk G$ with $\cV_Q^+$. Identify $M\bk K$ with $S^n$ as in \cite[Remark 2.8]{KelmerYu2023a}. We fix a measurable section $\mathfrak{s}:S^n\to K$ of the projection map from $K$ to $S^n$, and for any $\bm\alpha\in S^n$, we denote by $k_{\bm\alpha}=\mathfrak{s}(\bm\alpha)$ the corresponding representative.

Any $g\in G$ can be written as $g=u_{\bm x}a_yk$ with $k\in K$, and we denote its coordinates by $\bm x(g)$, $y(g)$ and $k(g)$. We normalize Haar measure by
\begin{align}\label{equ:Haar1}
\d\mu_G(g)=\frac{\d\bm x\,\d y}{y^{n+1}}\,\d\mu_K(k),
\end{align}
where $\d\bm x$ is Lebesgue measure on $\R^n$ and $\mu_K$ is Haar probability measure on $K$. The subgroup $L={U}M$ is unimodular, and we normalize its Haar measure by
\begin{align}\label{equ:Lhaar}
\d\mu_L(u_{\bm x}m)=\d\bm x\,\d\mu_M(m),
\end{align}
where $\mu_M$ is Haar probability measure on $M$.

Every element of $\cV$ can be written uniquely as $La_yk_{\bm\alpha}$ with $y>0$ and $\bm\alpha\in S^n$. The corresponding quotient measure is
\begin{align}\label{equ:lebmea}
\d m_{\cV}(a_yk_{\bm\alpha})=y^{-(n+1)}\,\d y\,\d\sigma_n(\bm\alpha),
\end{align}
where $\sigma_n$ is the rotation invariant probability measure on $S^n$. Under $Lg\mapsto\bm e_0g$, this is the cone measure $m_{\cV_Q^+}$ of \cite[(2.7)]{KelmerYu2023a}.

\subsection{Cusps of $\Gamma$}
 Let $P={U}AM$ denote the parabolic subgroup with ${U},A,M$ as above. Any other parabolic subgroup of $G$ is of the form $P'=kPk^{-1}$ for some $k\in K$ and its unipotent radical is given by ${U}'=k{U}k^{-1}$. The cusps of $\G$ are the $\G$-conjugacy classes of parabolic subgroups of $G$ whose unipotent radicals intersect $\G$ nontrivially. Let $\kappa$ denote the number of cusps of $\G$ and let $P_1,\ldots, P_{\kappa}$ be a full set of representatives of these classes. Note that for each $P_i$, we can take $k_i\in K$ such that $P_i=k_iPk_i^{-1}$, and thus it has a Langlands decomposition $P_i={U}_iA_iM_i$ with ${U}_i=k_i{U}k_i^{-1}$, $A_i=k_iAk_i^{-1}$ and $M_i=k_iMk_i^{-1}$. Write $\G_{{U}_i}:=\G\cap {U}_i$ and $\G_{P_i}:=\G\cap P_i$. By definition, $\G_{{U}_i}$ is non-trivial, and it is a finite-index subgroup of $\G_{P_i}$; see \cite[Lemma 2.1]{KelmerYu2023b}. In particular, since $A_i$ has no nontrivial finite subgroup, this implies that $\Gamma_{P_i}\subset L_i{:={U}_iM_i=k_iLk_i^{-1}}$.
For each $P_i$, we fix the scaling matrix $\tau_i=k_ia_{y_i}$ with $k_i\in K$ as above and $y_i>0$ the unique number such that $\mu_L(\tau_i^{-1}\G_{P_i}\tau_i\bk L)=1$.

\subsection{Eisenstein series and incomplete Eisenstein series}
The (spherical) Eisenstein series corresponding to the $i$-th cusp $P_i$ is defined,
for $\Re(s)>n$ and $g\in G$, by the convergent series
\begin{equation}\label{e:EisensteinSph}
	E_i(s, g)=\sum_{\gamma\in \Gamma_{\!P_i}\bk \Gamma} y(\tau_i^{-1}\gamma g)^s.
\end{equation}
The series $E_i(s, g)$ has a meromorphic continuation to the whole $s$-plane, which on the half plane $\Re(s)\geq \frac{n}{2}$ is holomorphic except for a simple pole at $s=n$ (called the trivial pole) with constant residue ${\omega_\Gamma:=\frac{1}{\mu_G(\Gamma\bk G)}}$, and possibly finitely many simple poles on the interval $(\frac{n}{2}, n)$ (called exceptional poles). We denote by $\cC_{\Gamma}\subseteq (\frac{n}{2}, n)$ the finite set of exceptional poles of Eisenstein series in all cusps. The constant term of $E_i(s,g)$ along the $j$th cusp is given by
$$E_{ij}(s,g)=\frac{1}{\vol(\R^n/\Lambda_j)}\int_{\R^n/\Lambda_j} E_i(s,\tau_ju_{\bm{x}}g)\,\d\bm{x},$$
where $\Lambda_j=\{\bm{x}\in \R^n: u_{\bm{x}}\in \tau_j^{-1}\Gamma_{{U}_j}\tau_j\}$. It satisfies
\begin{equation}
E_{ij}(s,g)=\delta_{ij}y(g)^s+\phi_{ij}(s)y(g)^{n-s}.
\end{equation}
{Here $\delta_{ij}$ is the Kronecker delta and $\Phi(s)=(\phi_{ij}(s))_{i,j=1}^\kappa$ is the scattering matrix for the chosen scaling matrices.}

For any bounded and compactly supported function $f: \cV\to \C$ that we may think of as a left $L$-invariant function on $G$, the incomplete Eisenstein series attached to $f$ at $P_i$ is defined by the series
\begin{align*}
E^{(i)}_f(g):=\sum_{\gamma\in \G_{P_i}\bk \G}f(\tau_i^{-1}\gamma g), \quad (g\in G),
\end{align*}
where $\G_{P_i}=\G\cap P_i$ and $\tau_i$ is the $i$-th scaling matrix as above.
This series is absolutely convergent and defines a function in $L^2(\Gamma\bk G)$.
\begin{remark}

{Suppose now that $Q$ is rational and $\Q$-isotropic, and let $\Gamma_Q:=\SO_Q(\Z)\cap G$. For $\Gamma=\Gamma_Q$, the light-cone Siegel transform $S_Qf$ from Section \ref{sec:recall} can be written as}
\begin{equation}\label{e:SQ-cusps}
S_Qf=\sum_{i=1}^\kappa E^{(i)}_{f_{\lambda_i}},
\end{equation}
where $f_{\lambda_i}(\bm{v})=f(\lambda_i^{-1} \bm{v})$ and each $\lambda_i>0$ is a constant related to the volume of the $i$th cusp; {see \cite[Lemma 2.1 and Remark 2.19]{KelmerYu2023a}}.
\end{remark}

\subsection{Inner products of incomplete Eisenstein series}\label{subsec:inner-products}
We can now record the formula obtained in \cite{YuThesis18} for the second moment of incomplete Eisenstein series, or rather its generalization given in \cite{KelmerYu2023a} for the inner product of two incomplete Eisenstein series corresponding to two cusps.

We first need to define the sesquilinear forms $M_{f,f'}(s)$ appearing as secondary terms in this formula.
For any $m\geq 0$ let $\cH_m\subseteq L^2(S^n,\sigma_n)$ denote the space of degree $m$ spherical harmonics and let $\{Y_{m,\ell}\}_{\ell=1}^{\dim\cH_m}$ be an orthonormal basis for $\cH_m$. For any bounded, compactly supported function $f:\cV\to \C$ let
\begin{equation}\label{e:fml}
 f_{m,\ell}(y)=\int_{S^n}f(a_y k_{\bm{\alpha}}){\overline{Y_{m,\ell}(\bm{\alpha})}}\,\d\sigma_n(\bm{\alpha}),
\end{equation}
and let $\hat{f}_{m,\ell}(s)$ denote their Mellin transforms
$$
\hat{f}_{m,{\ell}}(s)=\int_0^\infty f_{m,{\ell}}(y)y^{-s}\,\frac{\d y}{y}.
$$
For any $m\geq 0$ we define the functions  $P_m(s)$ by setting $P_0(s)=1$ and for $m\geq 1$
\begin{equation}\label{e: Pm}
P_m(s)=\prod_{j=0}^{m-1}\frac{n-s+j}{s+j}=\frac{\Gamma(n-s+m)\Gamma(s)}{\Gamma(n-s)\Gamma(s+m)}.
\end{equation}
{With this notation, for any $s\in (\frac{n}{2},n)$ and for any bounded, compactly supported functions $f,f':\cV\to\C$, we define}
\begin{equation}\label{e:Ms}
M_{f,f'}(s)=\sum_{m=0}^\infty \sum_{\ell=1}^{\dim(\cH_m)} P_m(s)\hat{f}_{m,\ell}(s)\overline{\hat{f}_{m,\ell}'(s)}.
\end{equation}

{Let $\mu_\Gamma$ denote the $G$-invariant probability measure on $\Gamma\backslash G$, normalized by $\d\mu_\Gamma=\omega_\Gamma\,\d\mu_G$.}
We can now state the result of \cite[Theorem 2.3]{KelmerYu2023a}, giving the following inner product formula {for $f_1,f_2\in C_c^\infty(\cV)$ and $1\leq i,j\leq\kappa$}:
\begin{align}\label{e:moment}
\langle E^{(i)}_{f_1},E^{(j)}_{f_2}\rangle_{L^2(\G\bk G,\mu_{\G})}
&=\omega_\Gamma^2 m_{\cV}(f_1)m_{\cV}(\overline{f_2})
+\omega_\Gamma\langle\delta_{ij}f_1+\cT_{ij}(f_1),f_2\rangle_{L^2(\cV,m_{\cV})}\nonumber\\
&+\omega_\Gamma\sum_{s_k\in\cC_\G}M_{f_1,f_2}(s_k)\Res_{s=s_k}\phi_{ij}(s).
\end{align}
Here each $\cT_{ij}:L^2(\cV,m_{\cV})\to L^2(\cV,m_{\cV})$ is a bounded linear operator {with operator norm $\|\cT_{ij}\|_{\mathrm{op}}\leq1$.}

\begin{remark}
{The moment formula \eqref{e:moment1}  for the light-cone Siegel transform follows from the inner product formula \eqref{e:moment}, the decomposition \eqref{e:SQ-cusps} and the fact that $\|\cT_{ij}\|_{\mathrm{op}}\leq1$ for each $1\leq i,j\leq \kappa$; see \cite[Section 2.4]{KelmerYu2023a} for details.}
\end{remark}

\section{Estimates for \texorpdfstring{$M_{f,f'}(s)$}{M(f,f';s)}}

\subsection{An integral formula}
In this section we prove \thmref{t:Msbound}.
The first step in the proof is to give an alternative formula for $M_{f,f}(s)$ that does not depend on the individual projections $f_{m,{\ell}}(y)$.
For any bounded compactly supported $f:\cV\to \C$ and any $s\in (\frac{n}{2},n)$ let  $H_sf\in L^2(S^n)$ be given by
\begin{equation}
H_sf(\bm{\alpha})=\int_0^\infty f(a_yk_{\bm{\alpha}})y^{-s}\,\frac{\d y}{y}.
\end{equation}
We then have the following integral formula for $M_{f,f'}(s)$.
\begin{Prop}
For any $s\in (\frac{n}{2},n)$, for any bounded compactly supported functions $f,f':\cV\to \C$
\begin{equation}\label{e:int}
M_{f,f'}(s)=\frac{1}{c(s,n)}\int_{S^n}\int_{S^n}\frac{H_sf(\bm{\alpha})\overline{H_sf'(\bm{\alpha}')}}{\|\bm{\alpha}'-\bm{\alpha}\|^{2(n-s)}}\,\mathrm{d}\sigma_n(\bm{\alpha})\,\mathrm{d}\sigma_n(\bm{\alpha}'),
\end{equation}
where $c(s,n)=2^{2(s-n)}\frac{\Gamma(n)\Gamma(s-\frac n2)}{\Gamma(\frac n2)\Gamma(s)}$.
\end{Prop}
\begin{proof}
{Throughout this proof, $\langle\cdot,\cdot\rangle$ denotes the inner product on $L^2(S^n,\sigma_n)$.
It suffices to prove the identity when $f=f'$; the general case follows by applying this case to $f\pm f'$ and $f\pm i f'$, and expanding.}

Next note that by expanding the integral and changing order of integration we see that
$$\langle H_sf, Y_{m,\ell}\rangle=\hat{f}_{m,\ell}(s),$$
hence
$$M_{f,f}(s)=\sum_{m=0}^\infty P_m(s)\sum_{\ell=1}^{\dim{\cH_m}}|\langle H_sf, Y_{m,\ell}\rangle|^2.$$
For any $m\geq 0$ let $\pi_m:L^2(S^n, \sigma_n)\to \cH_m$ denote the orthogonal projection to $\cH_m$ and using that $\langle H_sf, Y_{m,\ell}\rangle=\langle \pi_mH_sf, Y_{m,\ell}\rangle$ and that for any $h\in \cH_m$ we have $\sum_{\ell} |\langle h,Y_{m,\ell}\rangle|^2=\langle h,h\rangle$ we can further simplify this to
\begin{align*}
M_{f,f}(s)&=\sum_{m=0}^\infty P_m(s) \langle \pi_m H_sf, \pi_m H_s f\rangle
=\sum_{m=0}^\infty P_m(s) \langle \pi_m H_sf, H_s f\rangle.
\end{align*}

The next key idea is to identify $P_m(s)$ as an eigenvalue of an integral operator acting on $\cH_m$. Explicitly, by using the Funk-Hecke identity applied to the function $(1-t)^{-\lambda}$, it was shown in
\cite[Corollary 4.3]{FrankLieb2012} that for any {$0<\lambda<\frac{n}{2}$} and $h\in \cH_m$,
$$
\int_{S^n} \frac{h(\bm{\alpha})}{(1-\bm{\alpha}'\cdot\bm{\alpha})^\lambda}\,\d\sigma_n(\bm{\alpha})=2^{-\lambda}\frac{\Gamma(n)}{\Gamma(\frac{n}{2})}\frac{\Gamma(\frac{n}{2}-\lambda)\Gamma(m+\lambda)}{\Gamma(\lambda)\Gamma(n+m-\lambda)}h(\bm{\alpha}').
$$
Using this identity  with $\lambda=n-s$, recalling that $P_m(s)=\frac{\Gamma(s)}{ \Gamma(n-s)}\frac{\Gamma(m+n-s)}{\Gamma(m+s)}$, we get that
\begin{align*}
 \int_{S^n}\frac{h(\bm{\alpha})}{(1-\bm{\alpha}'\cdot\bm{\alpha})^{n-s}}\,\d \sigma_n(\bm{\alpha})={b(s,n)}P_m(s)h(\bm{\alpha}'),\quad (\forall\, h\in \cH_m),
\end{align*}
where ${b(s,n)}=2^{(s-n)}\frac{\Gamma(n) \Gamma(s-\frac{n}{2}) }{\Gamma(\frac{n}{2}) \Gamma(s) }$. {Thus $c(s,n)=2^{s-n}b(s,n)$.} Hence for any  $F\in L^2(S^n, \sigma_n)$  we have
$$P_m(s) \langle  \pi_m F,F \rangle=\frac{1}{{b(s,n)}}\int_{S^n}\int_{S^n}\frac{\pi_m F(\bm{\alpha})\overline{F(\bm{\alpha}')}}{(1-\bm{\alpha}'\cdot\bm{\alpha})^{n-s}}\,\d\sigma_n(\bm{\alpha})\,\d\sigma_n(\bm{\alpha}').$$
Summing over all $m$ we get that
$$\sum_{m=0}^\infty P_m(s) \langle  \pi_m F,F \rangle=\frac{1}{{b(s,n)}}\int_{S^n}\int_{S^n}\frac{F(\bm{\alpha})\overline{F(\bm{\alpha}')}}{(1-\bm{\alpha}'\cdot\bm{\alpha})^{n-s}}\,\d{\sigma_n}(\bm{\alpha})\,\d{\sigma_n}(\bm{\alpha}').$$

Applying this with $F=H_sf$ and noting that $\|\bm{\alpha}'-\bm{\alpha}\|^2=2(1-\bm{\alpha}'\cdot\bm{\alpha})$ concludes the proof.
\end{proof}
\subsection{Bounds for the integral formula}
We start with the following general integral bound.
\begin{Lem}\label{l:intbound}
For any bounded compactly supported function $h: {(0,\infty) \to [0,\infty)} $ and any $0<s<n$,
$$
\int_0^\infty h(y)y^{-s}\,\frac{\mathrm{d} y}{y}\leq  \frac{\|h\|_\infty^{1-\frac{s}{n}}n^{s/n}}{s}\left(\int_0^\infty h(y)y^{-n}\,\frac{\mathrm{d} y}{y}\right)^{s/n}.
$$
\end{Lem}
\begin{proof}
If $\|h\|_\infty=0$, both integrals vanish. Otherwise divide $h$ by $\|h\|_\infty$ and, changing it on a null set if necessary, we may assume $0\leq h\leq1$.
Making the change of variables $u=y^{-n}$ and setting $F(u)=h(u^{-1/n})$ (so that $0\leq F\leq 1$), we need to show that
 $$
 \frac{1}{n}\int_0^\infty F(u) u^{\frac{s}{n}-1}\,\d u\leq \frac{1}{s}\left(\int_0^\infty F(u)\,\d u\right)^{s/n}.
 $$
Since $u^{\frac{s}{n}-1}$ is decreasing and $0\leq F(u)\leq 1$, the bathtub principle
\cite[Theorem 1.14]{LiebLoss2001} shows that $\int_0^\infty F(u) u^{\frac{s}{n}-1}\,\d u$ is maximized when we replace $F$ by
$\chi_{[0,T]}$ with $T= \int_0^\infty F(u)\,\d u$, hence
$$
\frac{1}{n}\int_0^\infty F(u)\,u^{\frac{s}{n}-1}\,\d u\leq \frac{1}{n}\int_0^Tu^{\frac{s}{n}-1}\,\d u=\frac{1}{s}\Big(\int_0^\infty F(u)\,\d u\Big)^{s/n},
$$
as desired.
\end{proof}

We can now prove \thmref{t:Msbound}.
\begin{proof}
Let $\tfrac{n}{2}<s<n$ and let $p=\frac{n}{s}$.
For any bounded compactly supported function $f:\cV\to \C$ and  for any $\bm{\alpha}\in S^n$, {let $h_{\bm\alpha}(y):=|f(a_yk_{\bm{\alpha}})|$. Note that we have $\|h_{\bm\alpha}\|_{\infty}\leq \|f\|_{\infty}$ for $\sigma_n$-a.e. $\bm\alpha\in S^n$. For any $\bm\alpha\in S^n$ with $\|h_{\bm\alpha}\|_{\infty}\leq \|f\|_{\infty}$}, we apply \lemref{l:intbound} to $h_{\bm\alpha}$ to get that
$$
|H_sf(\bm{\alpha})|^p\leq \left(\int_0^\infty h_{\bm\alpha}(y)y^{-s}\,\frac{\d y}{y}\right)^p\ll_{n,s} \|f\|_\infty^{p-1}\int_0^\infty |f(a_yk_{\bm{\alpha}})| y^{-n}\,\frac{\d y}{y},$$
where the implied constant depends on $n,s$ but not on $\bm{\alpha}$.
Since the above inequality holds almost surely for $\bm\alpha\in S^n$, we can integrate it over this full measure subset of $S^n$ to get
\begin{align*}
\|H_sf\|_{L^p(S^n,\sigma_n)}^p
&\ll_{s,n} \|f\|_\infty^{p-1}\int_{S^n}\int_0^\infty |f(a_y k_{\bm{\alpha}})| y^{-n} \,\frac{\d y}{y}\,\d\sigma_n(\bm{\alpha})=\|f\|_\infty^{p-1} \|f\|_1.
\end{align*}

Finally, by \eqref{e:int} and the spherical Hardy–Littlewood–Sobolev inequality (see \cite[Theorem 2.2]{FrankLieb2012}) we get
\begin{align*}
\left|M_{f,f'}(s)\right|&= \left|\frac{1}{c(s,n)} \int_{S^n}\int_{S^n}\frac{H_sf(\bm{\alpha})\overline{H_sf'(\bm{\alpha}')}}{\|\bm{\alpha}'-\bm{\alpha}\|^{2(n-s)}}\,\d\sigma_n(\bm{\alpha})\,\d\sigma_n(\bm{\alpha}')\right|\\
&\ll_{s,n}   \|H_sf\|_{L^{p}(S^n,\sigma_n)}\|H_sf'\|_{L^{p}(S^n,\sigma_n)}\\
&\ll  (\|f\|_\infty \|f'\|_\infty)^{1-1/p} (\|f\|_1 \|f'\|_1)^{1/p},
\end{align*}
as desired.
\end{proof}

\section{Applications}

\subsection{Discrepancy bound for points on the light cone}
\begin{proof}[Proof of \thmref{t:Dis}]
First let $f$ be smooth and compactly supported. If $c_Q>0$, then $\frac{n}{2}<s_n<n$ and \thmref{t:Msbound} gives
\[
c_QM_{f,f}(s_n)\ll_Q\|f\|_\infty^{2-\beta_Q}\|f\|_1^{\beta_Q}.
\]
Combining this with \eqref{e:moment1} proves \eqref{e:momentbound}. If $c_Q=0$, the exceptional term is absent and the sharper $O_Q(\|f\|_2^2)$ remainder follows directly.

For bounded, compactly supported $f$, choose smooth $f_j$ with a common compact support, $\|f_j\|_\infty\leq\|f\|_\infty$, and $f_j\to f$ in $L^1$ and $L^2$. Applying \eqref{e:momentbound} to $f_j-f_k$ shows that $S_Qf_j$ is Cauchy in $L^2$. The first-moment formula gives
\[
\|S_Qf_j-S_Qf\|_{L^1(\G_Q\bk G, \mu_Q)}\leq\omega_Q\|f_j-f\|_1,
\]
so the $L^2$ limit is $S_Qf$. Passing to the limit proves \eqref{e:momentbound}, including its sharper form when $c_Q=0$.

For a relatively compact Borel set $B$, taking $f=\chi_B$ we get
\begin{equation}\label{e:moment3}
\int_{\Gamma_Q\bk G}|S_Q\chi_B(g)|^2\,\d\mu_Q(g)
=\omega_Q^2m_{\cV_Q^+}(B)^2
+O_Q\!\left(m_{\cV_Q^+}(B)+m_{\cV_Q^+}(B)^{\beta_Q}\right).
\end{equation}
For a general finite-measure Borel set $B$, we can take an increasing family of relatively compact Borel sets $\{B_j\}_{j\in \N}$ with $B_j\nearrow B$. Thus $S_Q\chi_{B_j}\nearrow S_Q\chi_B$; then by the monotone convergence theorem we can extend both the first-moment formula and \eqref{e:moment3} to $B$. Subtracting the square of the mean gives \eqref{e:DisMoment}.
\end{proof}
\subsection{Discrepancy bounds for \texorpdfstring{$\Gamma$}{\Gamma}-orbits}
Given any nonuniform lattice $\Gamma\leq G$ for each cusp $P_i$ of $\Gamma$ consider the orbit $\cL_i=L\tau_i^{-1}\Gamma\subseteq \cV$.
For any finite-measure Borel set $B\subseteq \cV$ the corresponding discrepancy function is $D(\cL_i g, B)=|\#(\cL_i g\cap B)-\omega_\Gamma m_{\cV}(B)|$.
We then have the following.
\begin{Thm}\label{t:Dis3}
There is $c>0$ depending on $\Gamma$ such that for any $1\leq i\leq \kappa$,  for any finite-measure Borel set $B\subseteq \cV$
\begin{align}\label{equ:gaorbbd}
\int_{\G\bk G}D(\cL_i g,B)^2d\mu_\G(g)\leq c(m_{\cV}(B)+m_{\cV}(B)^{\frac{2s_\Gamma}{n}}),
\end{align}
where $s_\Gamma\in (\frac{n}{2}, n)$ is the largest exceptional pole for $\Gamma$; when $\G$ has no exceptional poles, $s_{\Gamma}=\frac{n}{2}$.
\end{Thm}
\begin{proof}
Combining \eqref{e:moment} with \thmref{t:Msbound} and using the same approximation argument as in the proof of \thmref{t:Dis}, for any bounded, compactly supported $f_1,f_2:\cV\to\C$ and $1\leq i,j\leq\kappa$, we obtain
\begin{align}\label{e:incomplete-momentbound}
\langle E^{(i)}_{f_1},E^{(j)}_{f_2}\rangle_{L^2(\Gamma\bk G,\mu_\Gamma)}
&=\omega_\Gamma^2m_{\cV}(f_1)m_{\cV}(\overline{f_2})\nonumber\\
&\quad+O_\Gamma\!\left(\|f_1\|_2\|f_2\|_2+
\sum_{s\in\cC_\Gamma}(\|f_1\|_\infty\|f_2\|_\infty)^{1-s/n}
(\|f_1\|_1\|f_2\|_1)^{s/n}\right).
\end{align}
If $\cC_\Gamma=\emptyset$, the error term is simply $O_\Gamma(\|f_1\|_2\|f_2\|_2)$. Taking $j=i$ and $f_1=f_2=\chi_B$ gives
$$\int_{\G\bk G}\#(\cL_i g\cap B)^2\, \d\mu_\G(g)=\omega_\Gamma^2m_{\cV}(B)^2+O_\Gamma\!\left(m_\cV(B)+m_\cV(B)^{2s_\Gamma/n}\right).$$
Subtracting the square of the mean and recalling that
$$\int_{\G\bk G}\#(\cL_i g\cap B)\, \d\mu_\G(g)=\omega_\G m_{\cV}(B),$$
concludes the proof.
\end{proof}

 \subsection{Remaining applications}

The remaining results follow from the arguments in
\cite{KelmerYu2023a}. Recall that $\cA_{\beta_Q,c}$, defined
in \cite[Section 3]{KelmerYu2023a}, consists of the finite-measure
Borel subsets of $\cV_Q^+$ satisfying
\begin{align*}
\int_{\Gamma_Q\bk G}D(\cL_Qg,B)^2\,\mathrm{d}\mu_Q(g)
\leq c\bigl(\tilde{m}_{\cV_Q^+}(B)+\tilde{m}_{\cV_Q^+}(B)^{\beta_Q}\bigr).
\end{align*}
Here $\tilde{m}_{\cV_Q^+}=\omega_Q m_{\cV_Q^+}$.
Thus by \thmref{t:Dis}, especially \eqref{e:DisMoment}, we can choose $c>0$
depending only on $Q$ so that this family contains all
finite-measure Borel sets.

\smallskip

 We first prove \thmref{thm:linearformintro}. Recall that \cite[Theorem 3.9]{KelmerYu2023a} assumes that all finite-measure Borel sets belong to $\cA_{1,c}$ for some $c>0$. If we replace this assumption by the assumption that they belong to $\cA_{\beta,c}$ for some $1\leq\beta<2$, a slight and straightforward modification of the same argument gives the same conclusion with the range $0<\nu<\min\{\delta,(a-b)/(d+3)\}$ replaced by
$0<\nu<\min\left\{\delta,\frac{(a-b)(2-\beta)}{d+3}\right\}$.  Here we use the notation of that theorem, and $\delta>0$
is a constant for which the volume estimate in \cite[(3.20)]{KelmerYu2023a} holds.
Taking $\beta=\beta_Q$ and following the proofs of \cite[Theorems 4.3 and 4.4]{KelmerYu2023a} extends those results without the assumption $c_Q=0$.
In particular, this proves \thmref{thm:linearformintro}.

\smallskip

The proof of \thmref{thm:khinttypegen} follows the proof of \cite[Theorem 1.9]{KelmerYu2023a}. The only change is in the proof of \cite[Theorem 5.7]{KelmerYu2023a}, where we apply \cite[Theorem 3.7]{KelmerYu2023a} with $\beta=\beta_Q$
instead of $\beta=1$. This gives the modified exponent in the remainder.

\smallskip

Finally, for \thmref{t:Veech} we take $n=1$ and use the isomorphism of $G=\SO^+(2,1)$ with $\operatorname{PSL}_2(\R)$.
Consider the standard action of $\SL_2(\R)$ on row vectors in $\R^2$ by linear transformations and note that the set $S_M$ of distinct holonomy vectors of oriented saddle connections is a finite union of discrete orbits $\bm{v}\Gamma_M$ of its Veech group $\Gamma_M\leq\mathrm{SL}_2(\R)$; see \cite[Section 1]{BurrinFairchild2026}.
For each discrete orbit $\bm{v}\G_M$, the stabilizer of $\bm{v}$ in $\SL_2(\R)$ is conjugate to $U$, while
the stabilizer of the line $\R\bm v$ is a parabolic
subgroup conjugate to $\{\pm I_2\}UA$, corresponding
to a cusp of $\Gamma_M$.
Since reversing the orientation of a saddle connection negates its
holonomy vector, we have $S_M=-S_M$. Let
$\pi:\R^2\setminus\{\bm{0}\}\longrightarrow
(\R^2\setminus\{\bm{0}\})/\{\pm1\}$
be the quotient map, and define the quotient measure by $\nu(A)=\frac12|\pi^{-1}(A)|$. Note that in this case $\cV\cong (\R^2\setminus\{\bm{0}\})/\{\pm1\}$ and $\nu$ is a constant multiple of $m_{\cV}$.
Write
$\pi(S_M)=\bigsqcup_{j=1}^{r}\cO_j$
as a disjoint union of discrete orbits of the image of $\Gamma_M$ in
$\PSL_2(\R)$. For each orbit we let $c_j>0$ denote the density constant determined by the first moment identity
$$
\int_{\Gamma_M\backslash\SL_2(\R)}\sum_{\bm v\in \pi^{-1}(\cO_j)}f(\bm vg)\,\d\mu_M(g)
=c_j\int_{\R^2}f(\bm x)\,\mathrm{d}\bm x,
\qquad (\forall\, f\in C_c(\R^2)).
$$
Here $\mu_M$ is the probability invariant measure on $\G_M\bk \SL_2(\R)$. Let $c_M=\sum_{j=1}^{r}c_j$ which is the Siegel-Veech constant for $S_M$.

For any set $B\subseteq \R^2$ we have
\begin{align*}
\#(S_Mg\cap B)=\sum_{j=1}^{r}\left(
\#\bigl(\cO_jg\cap\pi(B)\bigr)
+\#\bigl(\cO_jg\cap\pi(B\cap(-B))\bigr)
\right).
\end{align*}
Note that the quotient measure satisfies
$$
\nu(\pi(B))+\nu(\pi(B\cap(-B)))=\frac12|B\cup(-B)|+\frac12|B\cap(-B)|=|B|.
$$
Moreover, both quotient sets form increasing families when $B$ does, {and after identifying $\cV$ with $(\R^2\setminus\{\bm{0}\})/\{\pm1\}$, each $\cO_j$ corresponds to a fixed dilation of the discrete orbit $\cL_i$ in \thmref{t:Dis3} for some $1\leq i\leq \kappa$. Hence after rescaling the counting sets}, we can apply \thmref{t:Dis3} and the argument for \thmref{t:Dis2} to get for each $1\leq j\leq r$, for any $\e>0$, for a.e. $g\in \SL_2(\R)$ and for all $B\in \cB$ with $|B|$ sufficiently large,
\begin{align}\label{equ:cuspcounting}
\#\bigl(\cO_jg\cap\pi(B)\bigr)
+\#\bigl(\cO_jg\cap\pi(B\cap(-B))\bigr)=c_j |B|+O_{M,\cB, g,\e}(|B|^{\sigma_M+\e}).
\end{align}
 {We also note that it is possible for the second family $\{\pi(B\cap (-B))\}_{B\in \cB}$ to have uniformly bounded $\nu$-measure. In this case, by the above first moment identity the corresponding counting function is $O_g(1)$ almost surely, and we still get \eqref{equ:cuspcounting} using $\nu(\pi(B))=|B|+O(1)$ and absorbing the $O_g(1)$ term into the error in \eqref{equ:cuspcounting}.}
Summing \eqref{equ:cuspcounting} over $j$ recalling that $c_M=\sum_j c_j$ we get the desired estimate \eqref{equ:sl2ct} for $\#(S_Mg\cap B)$.

\bibliographystyle{alpha}

\bibliography{DKbibliog}

\end{document}